\documentclass{amsart}
\usepackage{amssymb, graphicx}
\usepackage[mathscr]{euscript}

\newtheorem{thm}{Theorem}[section]
\newtheorem{cor}[thm]{Corollary}
\newtheorem{lem}[thm]{Lemma}

\newtheorem{defnn}[thm]{Definition}

\newtheorem{remarkk}[thm]{Remark}

\newtheorem{examplee}[thm]{Example}

\newcommand{\fig}[2] { \includegraphics[scale=#1]{#2}  }

\def\zz{{\bf Z}}
\def\co{\colon\thinspace}

\def\zz{{\bf Z}}
\def\qq{{\bf Q}}

\def\rr{{\bf R}}
\def\pp{{\bf P}}

\newcommand{\calg}{\mathcal{G}}
\newcommand{\calc}{\mathcal{C}}

\title{Knot Mutation:  4--Genus of  Knots and Algebraic Concordance}
\author{Se-Goo Kim}
\author{Charles Livingston}
\date{\today}
\address{Department of Mathematics, University of California, Santa Barbara,
CA 93106}

\address{Department of Mathematics, Indiana University, Bloomington, IN 47401}

\email{sekim@math.ucsb.edu}

\email{livingst@indianal.edu}


\begin{document}

\begin{abstract}   Kearton observed that mutation can change the concordance
class of a knot.  A close examination of his example reveals that  it is of
4--genus 1 and has a mutant of 4--genus 0. The first goal of this paper is to
construct examples to show that for any pair of nonnegative integers $m$ and $n$
there is a knot of 4--genus
$m$ with a mutant of 4--genus $n$.  

A second result of this paper is a crossing change formula for the algebraic
concordance class of a knot, which is then applied to prove the
invariance of the algebraic concordance class under mutation.  The paper concludes
with an application of   crossing change formulas  to give a short new proof of
Long's theorem that strongly positive amphicheiral knots are algebraically
slice.  
 \end{abstract}

\maketitle
 

\section{Introduction}  The main goal of this paper is to examine  the effect of
knot mutation on two concordance invariants of knots, the 4--ball genus and the
algebraic concordance class.  In the first case the extent to which mutation can
change the 4--genus is completely described.  In the second case it is shown that
the algebraic concordance class of a knot, as defined by Levine~\cite{lev69b}, is
invariant under mutation.  In the course of our work we develop crossing change
formulas for algebraic knot invariants. In the final sections of this paper we
apply such an approach to demonstrate that    Long's theorem that strongly positive
amphicheiral knots are algebraically slice   is an immediate corollary of the
Hartley-Kawauchi theorem that such knots have Alexander polynomials that are
squares.  Lastly, we show that the Hartley-Kawauchi theorem also follows from a
similar crossing change approach.
 \vskip.1in
\noindent{\em Mutation and Algebraic Concordance.} The construction of a mutant
$K^*$ of a knot
$K$ 
 consists of removing a 3--ball $B$ from $S^3$ that meets $K$   in two proper
arcs and gluing it back in   via an involution $\tau$ of its boundary $S$,
where $\tau$   is orientation preserving and leaves the set $S\cap K$
invariant. This  is among the most subtle    constructions of knot theory in
that it leaves a wide range of   knot  invariants unchanged~\cite{ada89,
kaw94, kaw96,   kir89, kir-kla90, mey-rub90, ron94, rub87, rub99}.  Most
relevant to the work here is the statement of~\cite{cop-lic99} that the
Tristram-Levine signatures,
$\sigma_\omega$, are invariant under mutation, since,  for $\omega$ a prime
power root of unity,  these provide the strongest classical bounds on the
4--genus~\cite{mur65, tri69}:
$ |\sigma_\omega(K) | /2 \le g_4(K)$.  In Sections~\ref{algcon}
and~\ref{algcon3} we give proofs of  the more general result:

\begin{thm}\label{algcon2} For any knot $K$, $\phi(K) = \phi(K^*)$, where
$\phi:\calc \to \calg$ is Levine's homomorphism from the knot concordance
group to the algebraic concordance group~\cite{lev69b}.
\end{thm}

\noindent The proof of Section~\ref{algcon} is entirely self-contained and
in addition gives a previously unnoticed crossing change formula for the
algebraic concordance class of a knot.  (As a side note, in Section~\ref{smaph}
we use this crossing change formula to give a quick derivation of a result of
Long that strongly positive amphicheiral knots are algebraically slice.)  In
Section~\ref{algcon3} an alternate proof of Theorem~\ref{algcon2} is
presented; this argument is somewhat briefer, but depends on the detailed analysis
of Seifert forms given in~\cite{cop-lic99}.

\vskip.1in 
\noindent{\em Mutation and the 4-Genus of a Knot.} The 4--genus of a knot,
$g_4(K)$, is the minimum genus of an embedded surface bounded by $K$ in   the
4--ball. This can be defined in either the smooth or topological locally flat
category; the results of this paper apply in either. It is   an especially
challenging invariant to compute; there remain low crossing number knots for
which it is uncomputed, though  the smooth category 
 has advanced considerably in recent years, most notably with the solution of
the Milnor conjecture which gives the 4--genus of torus
knots~\cite{kro-mor93}. 

Almost nothing has been known concerning the interplay between mutation and
the 4--genus.   Basically the only success in this realm consists of 
Kearton's observation~\cite{kea89}   that an example of~\cite{liv83} yields
an example  for which mutation changes the concordance class of a knot.    A
close examination of that example shows that it has 4--genus 1, but it has a
mutant  of 4--genus 0.    Further such examples have since been developed
in~\cite{kir-liv99, kir-liv01}. Our main result regarding the 4--genus is the
following.

\begin{thm}\label{mainthm} For every pair of nonnegative integers  $m$ and
$n$, there is a knot $K$ with mutant $K^*$ satisfying $g_4(K) = m$ and
$g_4(K^*) = n$.
\end{thm}

It should be noted that the original argument of~\cite{liv83} was based on a
paper of Gilmer ~\cite{gil83a}  in which it is now known   an error appears. 
To correct for that, the argument of~\cite{liv83} should be based on a
3--fold branched cover rather than the 2--fold cover.  The present work thus
serves to give the corrected argument for~\cite{liv83}. 

\vskip.1in
\noindent{\em Strongly Positive Amphicheiral Knots.}  A knot $K$ is called strongly
positive amphicheiral if when viewed as a knot in
$\rr^3$ it has a representative that is invariant under the map of $\rr^3$,
$\tau(x,y,z) = (-x,-y,-z)$.  We consider two theorems:

\begin{thm}[Long's Theorem \cite{lon84}] If $K$ is strongly positive amphicheiral,
then $K$ is algebraically slice.
\end{thm}

\begin{thm}[Hartley-Kawauchi Theorem \cite{har-kaw79}] If $K$ is a strongly
positive amphicheiral knot, then the Alexander polynomial $\Delta_K(t) = (F(t))^2$,
where $F$ is a symmetric polynomial.
\end{thm}

In Section~\ref{smaph} we use crossing change formulas developed earlier to prove 
that Long's theorem is an immediate corollary of the Hartley-Kawauchi result.  In
Section~\ref{smaph2} we use a   crossing change argument to give a
new proof of the Hartley-Kawauchi theorem.

\section{Background on Casson-Gordon invariants}
 A key tool in the proof of Theorem~\ref{mainthm}   is the main theorem
from~\cite{gil83b}  bounding  Casson-Gordon invariants in terms of the
4--genus of a knot.  Here is a simplified description of that result, based
on the statement of the theorem and later remarks in~\cite{gil83b}.

\begin{thm}[Gilmer's Theorem] Let $K$ be an algebraically slice knot
such that
$g_4(K) = g$ and let  $M_q$ be the $q$--fold branched cover of $S^3$ branched
over $K$ with $q$ a prime power. Let $\beta$ denote the linking form on
$H_1(M_q,\zz)$.  Then $
\beta$ can be written as a direct sum $\beta_1 \oplus \beta_2$ such that 1)
$\beta_1$ has a presentation of rank $2(q-1)g$ and  2) $\beta_2$ has a
metabolizer $D$ such that for any character $\chi$ of prime power order on
$H_1(M_q, \zz)$ given by linking with an element in $D$, one
has:$$|\sigma(K,\chi)| \le 2qg.$$
\end{thm}

Here $\sigma(K,\chi)$ is the Casson-Gordon invariant, originally denoted
$\sigma_1
\tau(K,\chi)$ in~\cite{cas-gor86, gil83b}.  We will need to know that $D$ can
be taken to be equivariant with respect to the deck transformation of $M_q$. 
Details  concerning this and other points will be given below, as they arise.

In our applications the group $H_1(M_q,\zz)$ will also be a vector space over
a finite field, in which case a metabolizer for $\beta_2$ will be
half-dimensional.  Hence:

\begin{cor} In the statement of Gilmer's Theorem, if
$H_1(M_q, \zz) \cong H_1(M_q,\zz_p)$, a $\zz_p$--vector space, then the
conclusion  of (1) can be restated   as 1) $\dim( \beta_1) \le 2(q-1)g$ and
in  (2) the metabolizer $D$ satisfies $ \dim (D)
\ge (\dim(H_1(M_q, \zz_p)) - 2(q-1)g) / 2$.
\end{cor}


\section{The   Building Blocks}

\begin{figure}[h]
 \centerline{  \fig{.7}{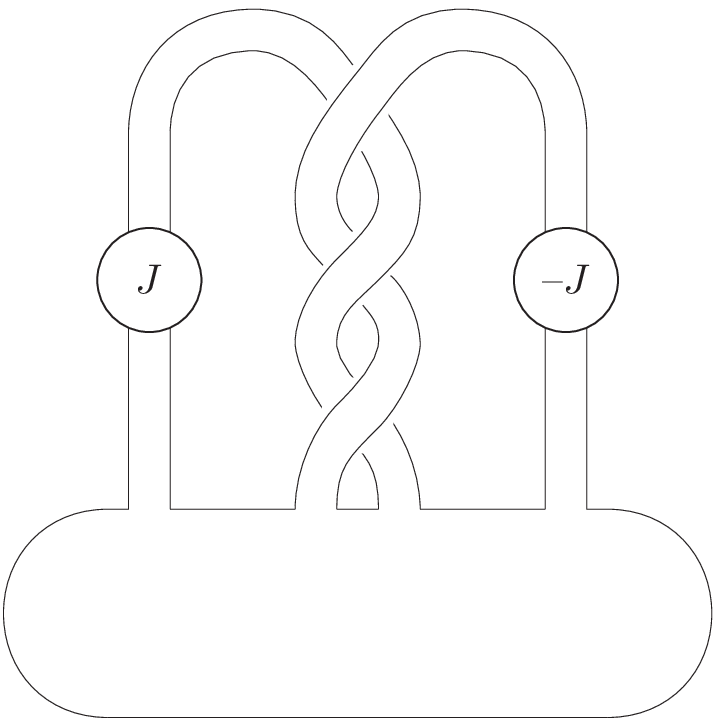}  }
 \vskip.05in
\caption{}\label{knot1} 
\end{figure}

Figure~\ref{knot1} illustrates a knot
$K_J$ of genus 1. The bands in the surface are tied in knots $J$ and
$-J$, for a knot $J$ to be determined later.  The twisting of the bands is
such that the Seifert matrix for 
$K_J$ is  
$    \left( \begin{matrix} 0 & 2 \\
               1 & 0 
\end{matrix} \right) . $  

Knots related to this one have been carefully analyzed elsewhere, for
example~\cite{gil-liv92, liv83,  liv01}, and the details of the following
results can be found there.  Here are the relevant facts.

\begin{enumerate}
\item  If $M_3$ denotes the 3--fold branched cover of $S^3$ branched over
$K_J$, then $H_1(M_3,\zz) = \zz_7 \oplus \zz_7$.

\item As a $\zz_7$--vector space, $H_1(M_3,\zz)$ splits as the direct sum of a
2--eigenspace, spanned by a vector $e_2$, and a 4--eigenspace, spanned by a
vector $e_4$,  with respect to the linear transformation induced by the deck
transformation.

\item  Linking with $e_i$ induces a character $\chi_i\co H_1(M_3,\zz) \to
\zz_7$.  Results of Litherland~\cite{lit84} (see also~\cite{gil93,
gil-liv92}) give the following: $$\sigma(K, \chi_2) = \sigma_{1/7}(J) +
\sigma_{2/7}(J) +
\sigma_{3/7}(J), $$ $$\sigma(K, \chi_4) = -\sigma_{1/7}(J) -
\sigma_{2/7}(J) -
\sigma_{3/7}(J). $$  Here $\sigma_{a/b}$ denotes the classical
Levine-Tristram signature, also written as  $\sigma_\omega$ with $\omega =
e^{(a/b)2 \pi i }$.  To simplify notation we abbreviate, for any knot $J$,
$$  s_7(J) = \sigma_{1/7}(J) + \sigma_{2/7}(J) +
\sigma_{3/7}(J). $$ There are knots for which $s_7$ is arbitrarily large, for
instance connected sums of trefoil knots. 
 
\end{enumerate}

\section{The Basic Examples}

We denote by $L_J$ the connected sum of $K_J $ with its mirror image,
reversed: $$L_J = K_J \# -K_J^r.$$ As observed by Kearton, $L_J$ is a mutant
of the slice knot $K_J \# -K_J$.

\begin{thm} For any choice of $J$, $g_4(L_J )\le 1$ and thus $g_4(n L_J) \le
n$. 
\end{thm}
\begin{proof} 
 
\begin{figure}[h]
 \centerline{  \fig{.65}{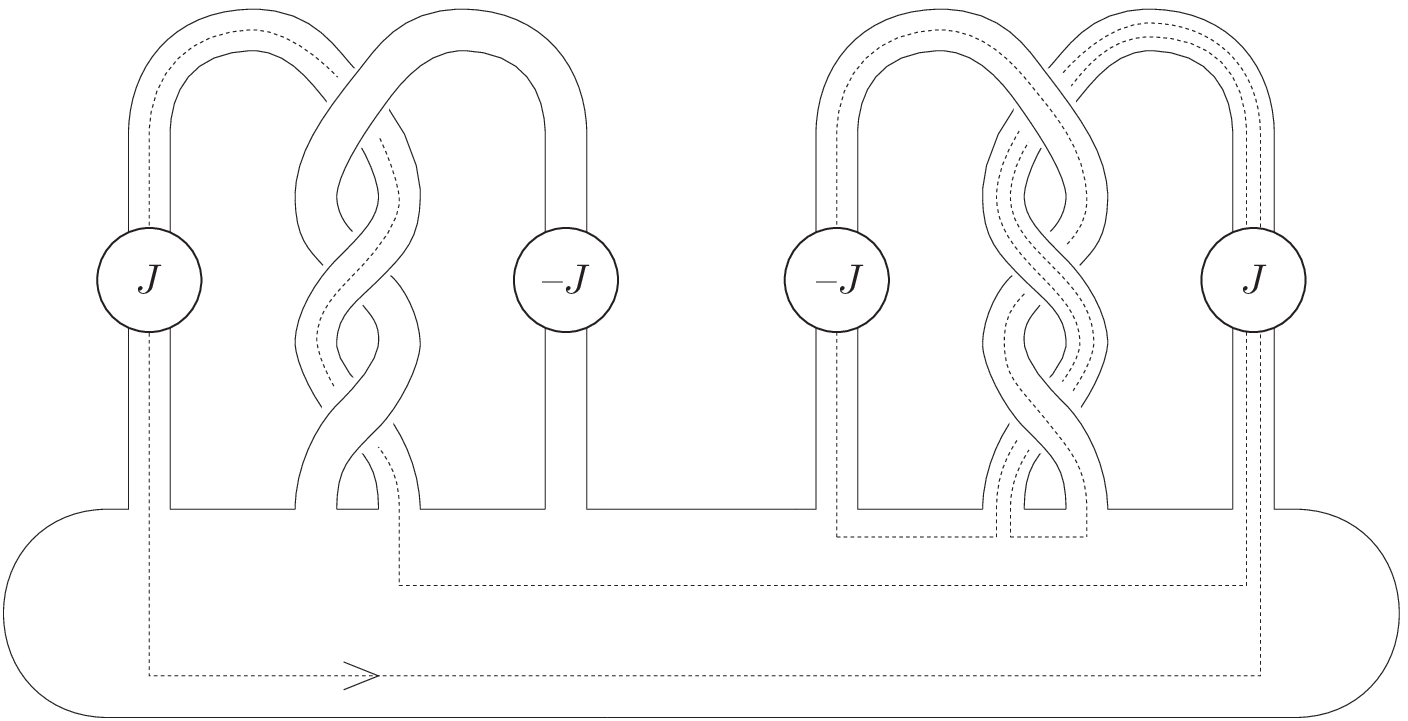}  }
 \vskip.05in
\caption{}  \label{knot2}
\end{figure}

Figure~\ref{knot2} illustrates $L_J$.  A simple closed curve on the genus 2
Seifert surface $F$ is indicated.  This curve has self-linking number 0 and
represents the slice knot $J  \# - J $.  Thus $F$ can be surgered in the
4--ball to reduce its genus to 1,  showing that $L_J$ bounds a surface of
genus 1 in the 4--ball, as desired.
\end{proof}

The homology of the 3--fold branched cover of $L_J$, $N_3$, naturally splits
as $(\zz_7 \oplus \zz_7) \oplus (\zz_7 \oplus \zz_7)$ with a 2--eigenspace
spanned by the vectors $e_2 \oplus 0 $ and $0 \oplus e_2' $, which we
abbreviate simply by
$e_2$ and $e_2'$.  Similarly for the 4--eigenspace.   We denote the
corresponding $\zz_7$--valued characters given by linking with $e_2$ and
$e_2'$ by $\chi_2$ and
$\chi_2'$, respectively. 

\begin{thm}The Casson-Gordon invariants of $L_J$ are given by:
$$\sigma(L_J , a \chi_2 + b \chi_2') = \epsilon(a)s_7(J) +
\epsilon(b)s_7(J) $$ 
$$\sigma(L_J , a \chi_4 + b \chi_4') = -\left( \epsilon(a)s_7(J) +
\epsilon(b)s_7(J)\right),$$ where $\epsilon(x) = 0$ or $1$ depending on
whether $x = 0$ or $x \ne 0$ mod 7.
\end{thm}

\begin{proof} This follows from the additivity of Casson-Gordon invariants;
see~\cite{lit84} or~\cite{gil83a}.  The only unexpected aspect of the formula
is that,  since we have $K_J \# - K_J^r$, it might have been anticipated that
the difference    $\epsilon(a)s_7(J) -
\epsilon(b)s_7(J)$ 
 would appear rather than the sum.  This switch occurs because the connected
sum involves the mirror image of reverse,
  rather than simply the mirror image; thus the role of $J$ and $-J$ are
reversed in the second summand.
\end{proof}


\section{Proof of  Theorem~\ref{mainthm}}

As observed by Kearton, for any knots $L_1$ and $L_2$, the connected sums $L_1
\# -L_2$ and $L_1 \#-L_2^r$ are mutants of each other.  Hence, it follows
immediately that for $m < n$, $n L_J$ is a mutant of $m L_J
\# (n-m)(K_J \# -K_J)$.  Since $ K_J \# - K_J$ is slice, this second knot is
concordant to, and hence of the same 4--genus as, $m L_J$.  To prove
Theorem~\ref{mainthm} we  show that for each positive integer
$n$ there exists a knot
$J$ so that for all $m \le n$, $g_4(m L_J) = m$.

At this point we fix some positive integer $n$ and select an arbitrary $m$, $1
\le m \le n$.  The knot  $J$  will be chosen as its necessary properties
become apparent.

Suppose that $mL_J$ bounds a surface $F$ in the 4-ball with genus
$g(F) = k < m$.  Let $V_3$ denote the 3--fold branched cover of $B^4$ branched
over $F$, having boundary the $m$--fold connected sum, $mN_3$.  Also,
abbreviate by
$D$ the image   of Tor($H_2(V_3, mN_3,
\zz)$) in
$H_1(mN_3, \zz)$.  An examination of the proof of Gilmer's Theorem
in~\cite{gil83b}  reveals that this $D$ is the metabolizer given in our
statement of his theorem above. Thus, for any $\chi$ corresponding to an
element in $D$ we have $|\sigma(mL_J, \chi)| \le 6k$.

With $\zz_7$--coefficients,  $H_1(mN_3, \zz)$ has dimension $4m$, so by
Gilmer's Theorem we have
$\dim(H_1(mN_3,
\zz) )- 2 {\dim}(D) \le 2(3-1)k = 4k$.  Hence, since
$k < m$, we have that $D$  is nontrivial.   

Observe that by its construction, $D$ is equivariant with respect to the deck
transformation   and hence contains an eigenvector.  Assume that it is a
2--eigenvector.  If we write
$H_1(mN_3, \zz) = \oplus_m H_1(N_3,\zz)$ then the 2-eigenvectors are
naturally denoted $e_{2,i}$ and $e_{2,i}'$, with $1 \le i \le m$, where
$e_{2,i}$ and $e_{2,i}'$ are the 2--eigenvector in the $i$th summand.  A
nontrivial 2--eigenvector in $D$ will be of the form
$\sum_i a_i e_{2,i} +
\sum_i b_i  e_{2,i}'$.  Using additivity, the Casson-Gordon invariant
corresponding to the dual character is given by:
$$ \left(\sum_i \epsilon(a_i)\right)  s_7(J) + \left(\sum_i \epsilon(b_i)
\right)s_7(J).$$  

To complete the proof, observe that this sum is greater than or equal to
$s_7(J)$, so that if $J$ is chosen so that $s_7(J) > 6n$ a contradiction is
achieved.  Notice that the choice of $J$ depends only on $n$ and not $m$.  

A similar argument applies if $D$ contains only a 4--eigenvector.


\section{The growth of $g_4(nK)$ for algebraically slice knots $K$.}

\begin{figure}[h]
 \centerline{  \fig{.6}{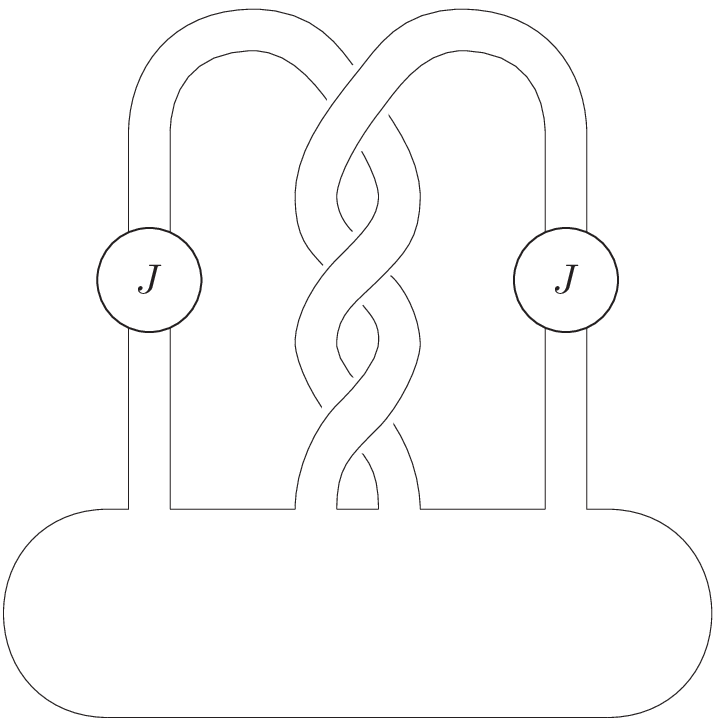}  }
 \vskip.05in
\caption{}\label{knot4} 
\end{figure}

For a  general knot $K$ one has $g_4(nK) \le n g_4(K)$ but one does not
usually have an equality.  In the case of a knot $T$, such as the trefoil,
for which the 4--genus is detected by a classical (additive) invariant, such
as the signature, one can sometimes demonstrate that
$g_4(nT) = ng_4(T)$.  But for algebraically slice knots with $g_4(K) \ne 0$
such arguments are not possible.  In fact, it is unknown whether in the
topological category there is such an algebraically slice knot
  for which the equality holds for all $n$.  (In the smooth setting the
second author, in~\cite{liv03}, has constructed an algebraically slice knot $K$ for which $g_4(K) = \tau(K) = 1$, where  $\tau$ is the invariant defined by  Ozv\'ath and
Szab\'o~\cite{osz-sza03}.  Since $\tau$ is  additive and bounds   $g_4$, it follows that $g_4(nK) = ng_4(K)$ for all $n$.)  We will here observe that one can come quite
close for the knot
$T_J$, where
$T_J$ is the knot illustrated in Figure~\ref{knot4}, built as
$K_J$ is, only with
$J$ tied in both bands rather than $J$ in one band and $-J$ in the other. 
(Similar results hold  for
$K_J$ and
$L_J$ but the proof would require the continued use of 3--fold covers rather
than the 2--fold cover for which the estimates are simpler.)

\begin{thm} For all $\epsilon$  with  $0 < \epsilon < 1$, there is a knot $J$
such that
$g_4(nT_J) > (1-\epsilon) n g_4(T_J)$ for all $n >0$.
\end{thm}

\begin{proof} Our proof builds upon Gilmer's original
argument~\cite{gil83b}.  Observe first that
$g_4(T_J)
\le 1$.  For the 2--fold branched cover we have that
$H_1(M_2, \zz) = \zz_3 \oplus \zz_3$ and the $\zz_3$--dimension  satisfies
$\dim (H_1(nM_2, \zz_3)) = 2n$.

If $nT_J$  bounds a surface in the 4--ball of genus $k$ less than or equal to
$(1-\epsilon)n$, then by Gilmer's theorem there exists a self--annihilating
summand $D$ with dim($H_1(nM_2, \zz_3)) - 2$dim$(D) \le 2k$ such that for all
characters $\chi$ dual to elements in $D$, one has $$|\sigma(nK_J, \chi)|
\le 4k.$$  

One computes that $\dim (D) \ge  n-k $.  A linear algebra  argument, basically
Gauss-Jordan elimination, now implies that some element of $D$ will be of the
form $\oplus_i \chi_i$ with at least $n-k$ of the $\chi_i$ nontrivial, and for
each of these $\chi_i$ the corresponding Casson-Gordon invariant is at least
$2\sigma_{1/3}(J)$.  Thus we have the equation
$$|(n-k)2\sigma_{1/3}(J)| \le 4k.$$

Since $k \le (1-\epsilon)n$, this reduces to 
$$|\epsilon n 2\sigma_{1/3}(J)| \le 4(1-\epsilon)n.$$ Simplifying yields $$| 
\sigma_{1/3}(J)| \le  2(1-\epsilon) /
\epsilon.$$  The proof is completed by noting that for any $\epsilon$ one can
select a $J$ for which this inequality does not hold.
\end{proof}


\section{Mutation and Algebraic Concordance}\label{algcon}

In this section we develop a crossing change formula for the
algebraic concordance class of a knot in order to prove
Theorem~\ref{algcon2}:  mutation preserves that algebraic concordance
class of a knot.     Certain knot invariants, such
as the Alexander polynomial and Tristram-Levine signatures, provide
algebraic concordance invariants, and these have been shown to be mutation
invariants (see for instance~\cite{ cop-lic99, lic-mil87}) but the
general question of whether mutation can change the algebraic concordance
class has remained open.  We should note that changing a knot to its
orientation reverse is a very  special case of mutation and reversal does
not change the algebraic concordance class of a knot, as   follows from
work of  Long~\cite{lon84}.  (More directly, it can be shown that the
complete set of algebraic concordance invariants defined by
Levine~\cite{lev69a} are unchanged by matrix transposition, the operation
on Seifert matrices induced by reversal.)

In the first subsection here we present a proof that the normalized Alexander
polynomial is invariant under mutation.  This argument is not new  but must
be presented to set up the needed notation for   the analysis of algebraic
concordance that follows.  The second subsection presents a review of the
theory and algebra of Levine's algebraic concordance group
$\calg$~\cite{lev69a}. In the final subsection we present a crossing change
formula for the algebraic concordance class of a knot and use this to prove
the mutation invariance of this class.

\vskip.1in

\subsection{The Alexander and Conway Polynomial} For an oriented  link $L$, a
choice of connected Seifert surface $F$ for $L$ and a choice of basis for
$H_1(F,\zz)$ there is a Seifert matrix $V(L)$, say of dimension $r \times
r$.  The (normalized) Alexander polynomial of $L$, $\Delta_L(t)$, can be
defined by setting
$$V_t(L)  = (1-t)V + (1-\bar{t})V^t  ,\ \  \mbox{and}$$  $$\Delta_L(t) =
\frac{1}{(z)^r}
\det(V_t(L)),$$ where $V^t$ denotes the transpose, $\bar{t} = t^{-1}$ and 
$z = t^{-1/2} - t^{1/2}$.  (Recall that $\Delta_L(t)$ can be expressed as a
polynomial in
$z$, $\Delta_L(t) = C_L(z) \in \zz[z]$, and this defines the Conway
polynomial~\cite{con70}.) Notice that
$z^2 = -(1-\bar{t})(1 -t)$, so that if $r$ is even (for instance, when $L$ is
connected, 
$r =  2\,  \mbox{genus}( F$))   $\Delta_L \in \zz[\bar{t},t]$ and elementary
  algebraic manipulations then  result in the usual normalized Alexander
polynomial: 
$$ \Delta_L(t) =t^{-r/2 } \det( V- tV^t).$$ (This polynomial is clearly
independent of change of basis and an observation below will show that it is
an
$S$--equivalence invariant~\cite{tro73} and thus depends only on $K$.)

 \begin{figure}[h]
 \centerline{  \fig{.5}{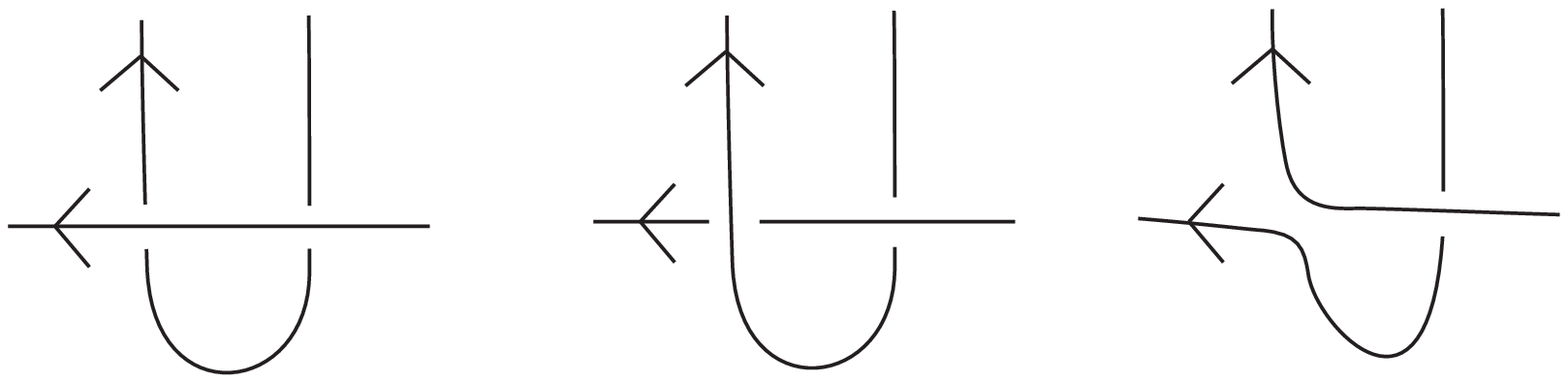}  }
 \vskip.05in
\centerline{  $L_-$ \hskip 1in  $L_+$  \hskip .95in $L_s$}
\caption{}\label{cchange} 
\end{figure} 

  Figure~\ref{cchange} illustrates a local picture of link diagrams for links
$L_-$, $L_+$, and
$L_s$, with the diagrams identical outside the local picture.  Any crossing
change and smoothing can be achieved using this local change.  In the diagram for $L_-$ a Reidemeister move eliminates the two crossings.  If Seifert's
algorithm is used to construct  a Seifert surface $F_0$  for $L_-$  using this simplified diagram, the corresponding    Seifert matrix   will be denoted $A$.  The Seifert surfaces for the links $L_-$ and $L_+$ that arise from Seifert's algorithm applied to the given diagrams are formed from $F_0$ by adding two twisted bands.  From this we have that
$V(L_\pm)$ is given by a $(r +2) \times (r+2)$ matrix of the form:

$$ V(L_{ \pm}) =
\left(
\begin{tabular}{ c   c}
 {\Huge{$A$}}& 

$   \begin{tabular}{c   c}                    
$a_1$ &0\\ 
$\vdots$ &$\vdots$\\
$a_r$ &0    \\ \end{tabular}
  $

 \\  

$   \begin{tabular}{c   c c}                    $a_1$  &$\cdots$ & $a_r$ \\ 
 0 & $\cdots$ & 0 \\ \end{tabular}
  $

 &   $   \begin{tabular}{c   c }                    $b$  &$1$ \\ 
 0 & $\epsilon_\pm$ \\ \end{tabular}
  $  \\
\end{tabular}
\right),
$$ where all entries are identical in these two matrices except that
$\epsilon_- = 0$  and
$\epsilon_+ = -1$.  Also, $V(L_s)$ is given by the same matrix, except with
the last row and column deleted.

A few consequences of these calculations follow quickly.

\begin{thm}The normalized Alexander polynomial is an $S$--equivalence
invariant and hence is a knot invariant.
\end{thm}
\begin{proof} $S$--equivalence is generated by the operation on Seifert
matrices which takes a matrix $A$ and replaces it with the matrix  denoted
$V(L_-)$ above.  That this doesn't change the Alexander polynomial is easily
checked: expand the relevant determinant along the last column and then along
the last row.
\end{proof}

\begin{thm}{\bf Conway skein relation.} The Alexander polynomial satisfies
$\Delta_{L_+} -
\Delta_{L_-} = z \Delta_{L_s}$.
\end{thm}

\begin{proof}  This again is a simple exercise in algebra, expanding the
determinant along the last column and then last row.
\end{proof}

\begin{thm}The Alexander polynomials of mutant knots are the same.
\end{thm}
\begin{proof}In the construction the mutant $K^*$, if the intersection of $K$
with the ball $B$ that is being taken out and replaced via an involution is
invariant under the extension of that involution to the 3--ball, then
$K^* = K$ and the polynomials are the same.  In general, a series of crossing
changes and smoothings converts $K \cap B$ into invariant tangles, so, via
the Conway skein relation, the polynomial of $K^*$ is the same as that for
$K$. 
\end{proof}

If $K$ is a knot, then the Alexander polynomial satisfies $\Delta_K(1) = 1$
and in  particular, $\Delta_K(t)$ is nontrivial.  Hence, in the above
matrices, working now with $K$ instead of $L$, $A_t$ is nonsingular.  Thus,
for $V_t(K_{\pm})$   the same set of row and column operations can be used to
eliminate the entries corresponding to the
$a_i$ in $V$.  There results  the
following  matrix, 
$W_t(K_\pm) $, where the entries are rational functions in
$t$, 
and   the matrix is hermitian with respect to the involution
induced by the map 
$t
\to \bar{t}$:
$$ W_t(K_{ \pm}) =
\left(
\begin{tabular}{ c   c}
 {\Huge{$A_t$}}& 

$   \begin{tabular}{ p{.5in}  p{.5in} }                    
 $0$ &\ \ 0\\ 
$\hskip.02in \vdots$ &$\hskip.11in \vdots$\\
$0$ &\ \  0    \\ \end{tabular}
  $

 \\  

$   \begin{tabular}{c   c c}                    $0$  &$\cdots$ & $0$ \\ 
 0 & $\cdots$ & 0 \\ \end{tabular}
  $

 &   $   \begin{tabular}{c   c }                    $c(t)$  &$1-t$ \\ 
 $1 -  \bar{t}$  & $\epsilon_\pm (1-t)(1- \bar{t})$ \\ \end{tabular}
  $  \\
\end{tabular}
\right),
$$

\begin{lem}\label{qupoly}     The ratio
$\Delta_{K_+} /
\Delta_{K_-}$  is equal to $c(t) +1$.
\end{lem}
\begin{proof}  This follows from a calculation of the relevant
determinants.\end{proof}

\vskip.2in

\subsection{Algebraic Concordance}

An algebraic Seifert matrix is a square integral matrix $V$ satisfying
$\det(V - V^t) = \pm 1$.  Such a matrix is called metabolic if it is
congruent to a matrix of the form 

$$\left(
 \begin{tabular}{c   c }                    
$0$  &$ A$ \\ 
 $B$  & $C$ \\ 

\end{tabular}
\right), 
$$  with $A, B,$ and $C$ square. Levine defined the algebraic concordance
group $\calg$ to be the set of equivalence classes of algebraic Seifert
matrices, with $V_1$ and $V_2$ equivalent if
$V_1 \oplus -V_2$ is metabolic.  The group operation is induced by direct sum.

A rational algebraic concordance group $\calg^\qq$ can be similarly defined,
where now it is
  required that $\det((V - V^t)(V + V^t)) \ne 0$.  Levine proved
in~\cite{lev69a} that the inclusion
$\calg
\to \calg^\qq$ is injective.

Consider next the set of nonsingular  hermitian  matrices with coefficients
in the field
$\qq(t)$, where  $\qq(t)$ has the involution $t \to \bar{t}$.  In this case
the equivalence relation generated by congruence to metabolic matrices
results in the Witt group of $\qq(t)$, denoted $W(\qq(t))$.

\begin{thm} The map $$V \to V_t = (1-t)V + (1 -
\bar{t} )V^t$$ induces an injection $ \calg \to W(\qq(t))$.
\end{thm}

\begin{proof} A proof is presented by Litherland~\cite{lit84} for
$\calg^\qq$, and the theorem follows from the injectivity of the inclusion
$\calg \to \calg^\qq$.  Note that in defining
$\calg^\qq$ (denoted $W_S(\qq, -)$ in~\cite{lit84})  Litherland restricts to
nonsingular matrices, but as he notes, Levine proved that every class in
$\calg$ has a nonsingular representative. To simplify notation, we will
use $W_t(K)$ to denote both the matrix and the Witt class represented by
the matrix when the meaning is clear in context.
\end{proof}

\subsection{Crossing Changes and Algebraic Concordance}

From the calculations and notation   above, if a crossing change is performed
on a knot $K$, the difference of Witt classes   associated to the Seifert
forms is given by 
$$W_t(K_+) - W_t(K_-) =  (A_t \oplus C_+) \oplus -(A_t \oplus C_-),$$ where
$$C_{\pm} =\left(
\begin{tabular}{c   c }                    $c(t)$  &$1-t$ \\ 
$1 -  \bar{t}$  & $\epsilon_\pm (1-t)(1- \bar{t})$ \\ \end{tabular}\right).$$

Since $A_t \oplus - A_t$ is Witt trivial, as is $C_-$, only $C_+$ contributes
to the difference of Witt classes.  Diagonalization, the identification of
$c(t) + 1$ with $\Delta_{L_+} / \Delta_{K_-}$, and a final multiplication of a
basis element (by $\Delta_{K_-}$) yields the following theorem.

\begin{thm}\label{cchange2} $W_t(K_+) - W_t(K_-)$ is represented by the matrix
$$  \left( \begin{tabular}{c   c }                    $\Delta_{K_+}(t)
\Delta_{K_-}(t)$  & $0$
\\ 
 $0$  & $-1$ \\ \end{tabular} \right),
  $$
  and thus the difference is determined by the Alexander polynomials of the
knots.
\end{thm}  

The special case of $\omega= -1$ in the following corollary is a result of
Murasugi~\cite{mur65}.  The proof follows from Theorem~\ref{cchange2} by
setting $t = \omega$ and induction on the number of crossing changes needed
to reduce $K$ to an unknot.  To avoid the matrix being nonsingular, we must
restrict to prime power roots of unity.

\begin{cor} For $\omega$ a prime power root of unity, $\mbox{\rm
sign}(\Delta_K(\omega) )= (-1)^{\sigma_\omega(K) /2}$.
\end{cor} We now have the main result of this section, the following
corollary of Theorem~\ref{cchange2},  a restatement of Theorem~\ref{algcon2}.

\begin{cor} The algebraic concordance class of a knot is invariant under
mutation; that is,
$W_t(K) = W_t(K^*)$ for any knot $K$ and its mutant $K^*$.  
\end{cor}

\begin{proof} A sequence of crossing changes in the tangle in $K$ that is
being mutated converts it into a tangle that is invariant under mutation. 
Thus we have a sequence of knots
$$K = K_0 , K_1, \ldots , K_n = K_n^* , K_{n-1}^*, \ldots , K_0^* = K^*,  $$
where $K_n = K_n^*$. By the previous theorem and the mutation invariance of
the Alexander polynomial, each pair of successive differences are equal: 
$W_t(K_i) - W_t(K_{i+1} )= W_t(K_i^*) - W_t(K_{i+1}^*)$.  Thus, $W_t(K) -
W_t(K_n) = W_t(K^*) - W_t(K_n^*)$.  Since $K_n = K_n^*$, the proof is
complete.
\end{proof}

\section{Generalized Mutation}\label{algcon3}

In~\cite{cop-lic99} Cooper and Lickorish study the effect of a generalization of
mutation, called {\it genus 2 mutation}, on the Seifert form of a knot.  Here we
deduce from their result an alternative proof of Theorem~\ref{algcon2}.  

Genus 2 mutation consists of removing a solid handlebody of genus 2 that
contains a knot $K$ from $S^3$ and replacing it via an involution of the
boundary. The involution is selected to extend to the solid handlebody so that
it has three fixed arcs.  The resulting knot is called $K^*$.  According
to~\cite{cop-lic99} there are Seifert matrices for $K$ and $K^*$ of the form

\[ V = 
\begin{pmatrix} A & B^t \\ B & C
\end{pmatrix}
\quad\text{and}\quad
V^\ast = 
\begin{pmatrix}
A & B^t \\ B & C^t
\end{pmatrix},
\]
respectively,
where $A$ and $C$ are square and $B$ is of the form $(0\mid b)$ for some
single column $b$.
Since $V$ is a Seifert matrix and $V-V^t=(A-A^t)\oplus (C-C^t)$,
$A$ and $C$ are also algebraic Seifert matrices.
Note that
\[
V_t=
\begin{pmatrix}
A_t & -z^2B^t \\ 
-z^2B & C_t
\end{pmatrix}
\quad\text{and}\quad
V^\ast_t=
\begin{pmatrix}
A_t & -z^2B^t \\ 
-z^2B & (C^t)_t
\end{pmatrix}
\]
where $z=t^{-1/2}-t^{1/2}$ and $z^2=-(1-t)(1-\bar{t})=-(1-t)-(1-\bar{t})$.

Since $A$ is a Seifert matrix, $A_t$ is nonsingular and hermitian.  Let
\[
P=
\begin{pmatrix}
I & z^2(A_t)^{-1}B^t \\ 0 & I
\end{pmatrix}.
\]
Then 
$V_t$ and $V^\ast_t$ are congruent to
$\bar{P}^t V_t P$ and $\bar{P}^t V^\ast_t P$, respectively,
which are, after a simple computation,
\[
\begin{pmatrix}
A_t & 0 \\
0 & C_t-z^4B(A_t)^{-1}B^t
\end{pmatrix}
\text{ and }
\begin{pmatrix}
A_t & 0 \\
0 & (C^t)_t-z^4B(A_t)^{-1}B^t
\end{pmatrix},
\]
respectively.
Suppose that $A$ is an $m\times m$ matrix.
Let $\alpha(t)\in\qq(t)$ be the $(m,m)$ entry of $(A_t)^{-1}$ and
recall that $B=(0\mid b)$ for some single column $b$ with
integral entries.
It is easy to see that
\[
B(A_t)^{-1}B^t = \alpha(t) bb^t.
\]
In particular, it is symmetric.  For simplicity, let
$E=C_t-z^4B(A_t)^{-1}B^t$.  Then 
$E^t = (C^t)_t - z^4B(A_t)^{-1}B^t$ and  we have that $V_t$ and
$V^\ast_t$ are congruent to
$A_t\oplus E$ and $A_t\oplus E^t$, respectively.
The difference of Witt classes of $V_t$ and $V^\ast_t$ is given by
\[
(A_t\oplus E)\oplus -(A_t\oplus E^t).
\]
Since $A_t\oplus -A_t$ is Witt trivial,
only $E\oplus -E^t$ contributes to the difference of Witt classes.
Observe that $E$ is a nonsingular hermitian matrix since $A_t\oplus E$ and
$A_t$ are. There is a nonsingular matrix $Q$ such that $F =
\bar{Q}^tEQ$ is diagonal.  This implies that $F=F^t=Q^t E^t\bar{Q}$. 
Using congruence by base change $Q\oplus \bar{Q}$, we see 
$E\oplus -E^t$ is congruent to
$F\oplus -F$, which is Witt trivial. 
Thus, $V_t=V^\ast_t$ in $W(\qq(t))$ and
$K$ and $K^\ast$ are algebraically concordant 
since $\calg \to W(\qq(t))$ is injective.

\section{Strongly Positive Amphicheiral Knots}\label{smaph}

A knot $K$ is called strongly positive amphicheiral if it is invariant
under an orientation reversing involution of $S^3$ that preserves the
orientation of $K$.  This is easily seen to be equivalent to the
statement that $K$, when viewed as a knot in
$\rr^3
\subset S^3$,   is isotopic to a knot, again denoted
$K$, which is invariant under the involution $\tau: \rr^3 \to \rr^3$
given by
$\tau(x) = -x$. 

Hartley and Kawauchi~\cite{har-kaw79} proved that if $K$ is strongly
positive amphicheiral then $\Delta_K(t) = (F(t))^2$ for some Alexander
polynomial $F$.  As a complementary result Long~\cite{lon84} proved that
strongly positive amphicheiral knots are algebraically slice.  Here we
demonstrate that Long's theorem is in fact a corollary of the
Hartley-Kawauchi theorem and the crossing change formula for the
algebraic concordance class.

A  bit of notation will be helpful: for a strongly amphicheiral knot that
is invariant under the involution $\tau$, $\tau$ defines a pairing of the
crossing points in a diagram of $K$.  A {\em paired crossing change} on
such a $K$ consists of changing both of a pair of crossings.  Notice that
since $\tau$ is orientation reversing, the two crossings will be of 
opposite sign, so we denote the original knot $K_{+ -}$ and the knot formed by
making the paired crossing changes
$K_{- +}$.

\begin{lem}  A sequence of paired crossing changes converts a strongly
positive amphicheiral knot into the unknot.

\end{lem}

\begin{proof} Since an involution of $S^1$ cannot have one fixed point,
$K$ misses the origin in $\rr^3$ and thus projects to a knot $\bar{K}$ in
the quotient $\rr^3 - \{0\} / \tau \equiv \rr \pp ^2 \times \rr$.  Since
$\bar{K}$ lifts to a single component in the cover, it is homotopic to standard
generator of $\pi_1 (\rr \pp ^2 \times \rr)$, whose lift is an unknot in
the cover.  That homotopy can be carried out by a sequence of crossing
changes, each of which lifts to a pair of crossing changes in the cover.
\end{proof}

\begin{thm}[Long's Theorem]  If $K$ is strongly positive amphicheiral,
then $K$ is algebraically slice.
\end{thm}

\begin{proof} By the previous lemma we need only show that $W_t(K_{+-}) -
W_t(K_{-+})$ represents 0 in $W(\qq(t))$.

Working in the Witt group we can write $$W_t(K_{+-}) - W_t(K_{-+}) =
(W_t(K_{+-}) - W_t(K_{--})) - (W_t(K_{-+}) - W_t(K_{--})).$$ 
Applying Theorem~\ref{cchange2}, this is represented by the difference
$$  \left( \begin{tabular}{c   c }                    $\Delta_{K_{+-}}(t)
\Delta_{K_{--}}(t)$  & $0$
\\ 
 $0$  & $-1$ \\ \end{tabular} \right) 
 \oplus 
-  \left( \begin{tabular}{c   c }                    $\Delta_{K_{-+}}(t)
\Delta_{K_{--}}(t)$  & $0$
\\ 
 $0$  & $-1$ \\ \end{tabular} \right)$$
Applying the Hartley-Kawauchi theorem, we write $\Delta_{K_{+-}}(t) = F(t)^2$ and
$\Delta_{K_{-+}}(t) = G(t)^2$, and then cancel  the $(-1)$ summands  to arrive
at  the difference
 $$  \left( \begin{tabular}{c   c }                    $F(t)^2
 \Delta_{K_{--}}(t)$  & $0$
 \\ 
  $0$  &   $-G(t)^2
 \Delta_{K_{--}}(t)$  \\ \end{tabular} \right). $$
This form has a metabolizer generated by the vector $(G(t), F(t) ) \in
\qq(t)^2$, and hence it is trivial in the Witt group, as desired.
\end{proof}

\section{The Hartley-Kawauchi Theorem}\label{smaph2}

Here we present a combinatorial proof of the theorem that for strongly
positive amphicheiral knots the Alexander polynomial is a square of an
Alexander polynomial.  The proof also gives an alternative, though
longer, route   to Long's theorem than was given in the previous section. 
We begin by considering the existence of an equivariant Seifert surface for
such a knot.

If Seifert's algorithm for constructing a Seifert surface is applied to a diagram
for a strongly amphicheiral knot that is invariant under $\tau$, the resulting
surface will be invariant. In addition, $\tau$ restricted to this surface
is orientation preserving since $\tau$ preserves the orientation of
the knot that is the boundary of the surface.  However $\tau$ reverses the positive normal
direction since it reverses the orientation of $R^3$.  Thus we have the following.

\begin{lem}\label{theta-tau} Let $K$ be a strongly positive amphicheiral knot with
involution
$\tau$.  Then a Seifert surface $F$ of $K$ can be constructed so that $F$ is
invariant under $\tau$ and its Seifert form $\theta$ satisfies $\theta(\tau
u, \tau v) = -\theta(v, u)$ for all $u , v \in H_1(F)$.

\end{lem}

\begin{figure}[h]
 \centerline{  \fig{.5}{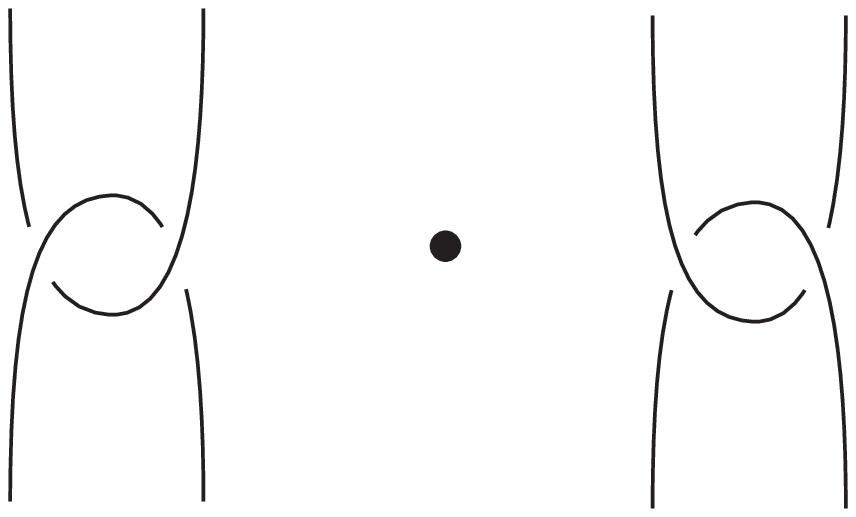}  }
 \vskip.05in
\caption{}\label{knot6} 
\end{figure}

To understand the effect of crossing changes, we consider
two figures.   Figure~\ref{knot6} represents a portion of a symmetric
diagram of a strongly amphicheiral knot, say
$K_{+-}$.  The dot in center of the figure represents the origin in $R^3$,
the center of symmetry.

\begin{figure}[h]
 \centerline{  \fig{.5}{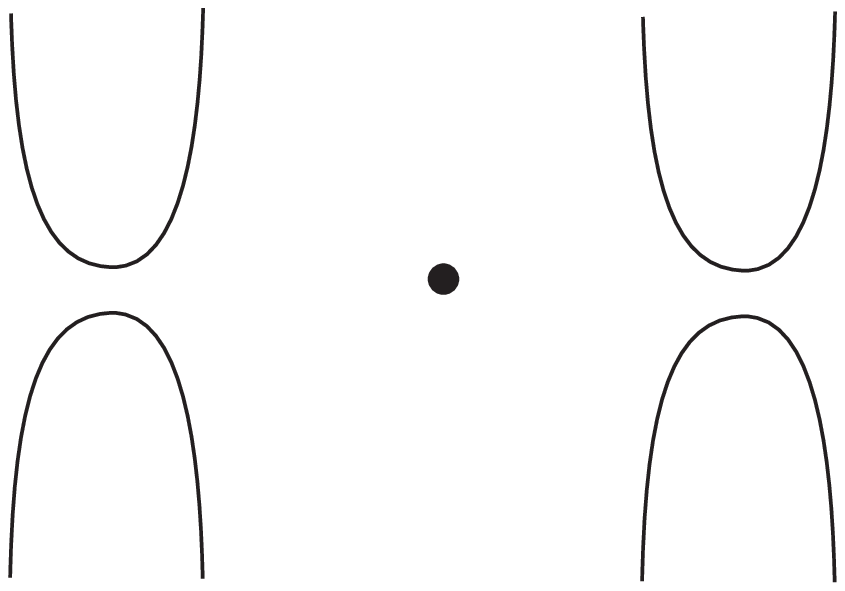}  }
 \vskip.05in
\caption{}\label{knot5} 
\end{figure}

For the  knot $K_{-+}$ the diagram will be the same, only a symmetric pair of
crossing changes has been made.  Thus, for $K_{-+}$ the clasps pull apart,
leaving a knot, denoted $K'$, with diagram as illustrated in
Figure~\ref{knot5}.

Suppose that $K'$ has an equivariant Seifert surface $F_0$ given by Seifert's
algorithm and $H_1(F_0)$ has symplectic basis $w_1, \ldots , w_r$.  Then an
equivariant Seifert surface $F$ for $K_{+-}$ is given by adding four bands to
$F_0$.  The basis for $H_1(F_0)$ can be naturally extended to symplectic one for
$H_1(F)$,
$w_1,
\ldots , w_r, x, y, \tau x, \tau y$, where $y$ has trivial Seifert pairing with
all elements other than $x$ and itself, and $x$ has trivial Seifert pairing  with
$\tau y$.

Let $A$ be the Seifert matrix of $F_0$ with respect to $w_1,\ldots,w_r$ and let
$T$ denote the matrix representing the action of $\tau$ on $H_1(F_0)$. Then
Lemma~\ref{theta-tau} applied to
$F_0$ can be rewritten in terms of matrices: $T^tAT=-A^t$.  After hermitianizing
and taking inverses, we have
\[
T(A_t)^{-1}T^t= -(A^t_t)^{-1}=\overline{-(A_t)^{-1}}.
\]

To find the Seifert matrix for $F$ with respect to the above basis, 
a couple of things have to be clarified.
First, note that
$\theta(x,\tau x)=-\theta(\tau\tau x,\tau x)=-\theta(x,\tau x)$ and
hence $\theta(x,\tau x)=0$.
Similarly, $\theta(\tau x,x)=0$.

Second, let $a =
\left(\begin{smallmatrix}
\theta(w_1,x) \\ \vdots \\ \theta(w_r,x)
\end{smallmatrix}\right)$ and
$T=(t_{ij})_{1\le i,j\le r}$.
Then
\[
\begin{pmatrix}
\theta(w_1,\tau x) \\ \vdots \\ \theta(w_r,\tau x)
\end{pmatrix}
=
\begin{pmatrix}
-\theta(x,\tau w_1) \\ \vdots \\ -\theta(x,\tau w_r)
\end{pmatrix}
=
\begin{pmatrix}
-\sum_j t_{j1} \theta(x,w_j) \\ \vdots \\ 
-\sum_j t_{jr}\theta(x, w_j)
\end{pmatrix}
=
-T^t
\begin{pmatrix}
\theta(x,w_1) \\ \vdots \\ \theta(x, w_r)
\end{pmatrix}
= -T^ta.
\]

It follows readily  that the Seifert matrix for $K_{+-}$ is the
$(r+4)\times (r+4)$ matrix:
\[
V^\epsilon =
\begin{pmatrix}
A & a & 0 & -T^ta & 0 \\
a^t & b & 1 & 0 & 0 \\
0 & 0 & \epsilon & 0 & 0 \\
-a^t T & 0 & 0 & -b & 0 \\
0 & 0 & 0 & -1 & -\epsilon
\end{pmatrix},
\]
where $\epsilon=-1$. 

Similarly, for $K_{-+}$ the same matrix arise, only in this case   $\epsilon=0$.
After hermitianizing
\[
V^\epsilon_t =
\begin{pmatrix}
A_t & -z^2a & 0 & z^2T^ta & 0 \\
-z^2a^t & -z^2b & 1-t & 0 & 0 \\
0 & 1-\bar{t} & -z^2\epsilon & 0 & 0 \\
z^2a^t T & 0 & 0 & z^2b & -(1-\bar{t}) \\
0 & 0 & 0 & -(1-t) & z^2\epsilon
\end{pmatrix},
\]
where $z=t^{-1/2}-t^{1/2}$.
Let
\[
P = 
\begin{pmatrix}
I & z^2(A_t)^{-1}a & 0 & -z^2(A_t)^{-1}T^ta & 0 \\
0 & 1 & 0 & 0 & 0 \\
0 & 0 & 1 & 0 & 0 \\
0 & 0 & 0 & 1 & 0 \\
0 & 0 & 0 & 0 & 1
\end{pmatrix}.
\]
Let $W^\epsilon_t = \bar{P}^t
V^\epsilon_t P$. Then
\[
W^\epsilon_t = 
\begin{pmatrix}
A_t & 0 & 0 & 0 & 0 \\
0 & -z^2b-z^4a^t(A_t)^{-1}a & 1-t & 
z^4a^t(A_t)^{-1}T^ta & 0 \\
0 & 1-\bar{t} & -z^2\epsilon & 0 & 0 \\
0 & z^4a^tT(A_t)^{-1}a & 0 & 
z^2b-z^4a^tT(A_t)^{-1}T^ta & -(1-\bar{t}) \\
0 & 0 & 0 & -(1-t) & z^2\epsilon
\end{pmatrix}.
\]
Let $c(t)=-z^2b-z^4a^t(A_t)^{-1}a$.
Since $W^\epsilon_t$ is hermitian, $c(t)=\overline{c(t)}$.
The (1,1) entry of the lower right $2\times 2$ submatrix of
$W^\epsilon_t$ is
\[
z^2b-z^4a^t\left(T(A_t)^{-1}T^t\right)a =
\overline{z^2b+z^4a^t(A_t)^{-1}a} =\overline{-c(t)}=-c(t).
\]
Let $d(t)=z^4a^t(A_t)^{-1}T^ta$.
Then the $1\times 1$ matrix $d(t)$ is equal to its transpose
\[
z^4a^tT(A^t_t)^{-1}a
=z^4a^tT\left(-T(A_t)^{-1}T^t\right)a
=-z^4a^t(A_t)^{-1}T^ta = -d(t)
\]
and hence $d(t)=0$. Also, note that
$z^4a^tT(A_t)^{-1}a = \overline{d(t)}=0$ since $W_t^\epsilon$ is
hermitian.

Thus $V^\epsilon_t$ is congruent to, by base change $P$,
\[
A_t\oplus C\oplus -C^t,
\text{ where }
C = \begin{pmatrix}
c(t) & 1-t \\
1-\bar{t} & -z^2\epsilon
\end{pmatrix}.
\]
Since $\det(P)=1$,
\[
\Delta_{K_{+-}}=(c(t)+1)^2\frac{1}{z^r} \det (A_t)
=(c(t)+1)^2 \Delta_{K_{-+}},
\]
where $c(t)=c(\bar{t})$.  This proves Hartley-Kawauchi
theorem.

Next, to prove Long's theorem, we will show that
$V_t(K_{+-})$, $A_t$, and
$V_t(K_{-+})$ are all Witt-equivalent.
It suffices to show that $C\oplus -C^t$ is Witt-trivial.
Observe that $C$ is nonsingular and hermitian since 
$A_t\oplus C\oplus -C^t$ and $A_t$ are.
There is a nonsingular matrix $Q$ such that $D=\bar{Q}^tCQ$ is diagonal.
This implies that $D=D^t=Q^t C^t\bar{Q}$. 
Using congruence by base change $Q\oplus \bar{Q}$, we see 
$C\oplus -C^t$ is congruent to
$D\oplus -D$, which is Witt trivial. 
Thus, $K_{+-}$ and
$K_{-+}$ are algebraically concordant.  This proves Long's theorem.


\end{document}